\RequirePackage{fix-cm}

\documentclass[preprint,12pt]{elsarticle}

\usepackage{epsfig}
\usepackage{fancybox}
\usepackage{multicol}
\usepackage{color}
\usepackage{graphicx}
\usepackage{url}
\usepackage{floatrow}
\usepackage{subfig}

\usepackage{amsmath}
\usepackage{amssymb}

\usepackage{amsthm}
\usepackage{mathrsfs} 
\usepackage{empheq}

\usepackage{lineno}

\makeatletter
\newcommand{\Rmnum}[1]{\expandafter\@slowromancap\romannumeral #1@}
\makeatother

\usepackage{multirow}

\usepackage{comment}
\usepackage[]{algorithm2e}

\begin{document}

\begin{frontmatter}

\title{High-order covariant differentiation in applications to Helmholtz-Hodge decomposition on curved surfaces}

\author{Sehun Chun}
\ead{sehun.chun@yonsei.ac.kr}
\address{Underwood International College, Yonsei University, South Korea}
\begin{abstract}
 A novel high-order numerical scheme is proposed to compute the covariant derivative, particularly for divergence and curl, on any curved surface. The proposed scheme does not require the construction of a curved axis or metric tensor, which would deteriorate the accuracy of the covariant derivative and prevent its application to complex surfaces. As an application, the Helmholtz-Hodge decomposition (HHD) is adapted in the context of the Galerkin method for displaying the irrotational, incompressible, and harmonic components of vectors on curved surfaces. 
\end{abstract}

\begin{keyword}
Covariant derivative \sep Moving frames \sep Connection form \sep Helmholtz-Hodge decomposition \sep Discontinuous Galerkin method
\end{keyword}

\end{frontmatter}

\section{Introduction} 

The covariant derivative on a curved surface is obtained differently from the Euclidean derivative because the axis is not fixed on the surface; moreover, the relative rotation of the axis is inevitable for the differentiation of a vector. For a vector $\mathbf{u}$, the covariant differentiation along a curved axis $x^{\alpha}$ is given as \cite{Moore}
\begin{equation}
\mathbf{u}_{;\alpha} = \frac{\partial \mathbf{u}}{\partial x^{\alpha}} =  \sum_{\mu} \left[ \frac{\partial u^{\mu}}{\partial x^{\alpha}} + \sum_{\nu} \Gamma^{\mu}_{\alpha \nu} u^{\nu} \right ] \boldsymbol{\nu}^{\mu},  \label{Covariant1}
\end{equation}
where $\boldsymbol{\nu}^{\mu}$ is the unit tangent vector of the axis $x^{\mu}$. The subscript $;\alpha$ indicates that the corresponding quantity is the covariant derivative with respect to the curved axis of the corresponding index $\alpha$. The variable $ \Gamma^{\mu}_{\alpha \nu}$ is referred to as the \textit{second type of Christoffel symbol} to represent how the axis $x^{\alpha}$ rotates as it moves along the curved axis $x^{\nu}$.

The first challenge in computing Eq. \eqref{Covariant1} for the general surface is to find a continuous and differentiable curved axis $x^{\alpha}$. This is especially difficult, both computationally and analytically, in regions with various curvatures or anisotropic properties caused by geometric singularity. The second challenge is to compute the corresponding metric tensor $g_{\alpha \beta}$ and the Christoffel symbol $\Gamma^{\mu}_{\alpha \nu}$ with a sufficiently small error in comparison to the discretization error. The Christoffel symbol can be directly obtained by differentiating the metric tensor $g_{\alpha \beta}$ as follows
\begin{equation*}
\Gamma^{\alpha}_{\mu \nu} = \frac{g^{\alpha \sigma}}{2}  \left [ \frac{\partial g_{\nu \alpha}}{\partial x_{\mu}} +  \frac{\partial g_{\sigma \mu}}{\partial x_{\nu}} -  \frac{\partial g_{\nu \mu}}{\partial x_{\sigma}} \right ],
\end{equation*}
where the tensor $g^{\alpha \beta}$ is the inverse of $g_{\alpha \beta}$. Computing the metric tensor is particularly challenging because it requires computing the length of the curved axis. Invalid construction of the curved axis or inaccurate computation of the length of the axis yields Christoffel symbols with non-negligible errors. Inaccurate metric tensor and Christoffel symbol function as corrupted coefficients of partial differential equations (PDEs) that cause nonphysical dynamics. 

In image processing and surface PDEs, the covariant derivative has been a crucial tool. Recent works on the computation and application of covariant derivatives are as follows: the computation of covariant derivative by discrete connection on triangulated 2-manifold \cite{Liu}, application of covariant derivative to image regularization \cite{Batard}, and the reformation of covariant derivative in Cartesian coordinates in the context of finite element methods \cite{Nestler}. Extensive literature of applications and comparisons to covariant differentiation for diffusion equations and the shallow water equations can be found in ref. \cite{MMF2} and \cite{MMF3}, respectively.

This paper introduces a novel method of computing a high-order covariant differentiation in Eq. \eqref{Covariant1} without constructing a curved axis $x_{\alpha}$ or a Christoffel symbol  $\Gamma^{\mu}_{\alpha \nu}$. To achieve this, we introduce moving frames and their special arrangement, known as the \textit{connection form}.

\section{Connection form}
For $1 \le i \le 3$, let $\mathbf{e}^i$ be moving frames constructed at each point to constitute a tetrahedron. Frames $\mathbf{e}^i$ are orthonormal such that $\mathbf{e}^i \cdot \mathbf{e}^j$ = $\delta^i_j$ where $\delta^i_j$ is the Kronecker delta. Let $\Omega_i$ be the tessellation of a smooth surface $\Omega$ such that $\cup \Omega_i = \Omega$ and $\Omega_i \cap \Omega_j = \delta^i_j$. Let $\Omega_i$ be locally Euclidean such that an orthogonal axis can be built at every point. $\mathbf{e}^i$ is differentiable in each element $\Omega_i$ but may not be differentiable across the interfaces. For constructions and more details on moving frames, refer to refs. \cite{MMF1, MMF2, MMF3, MMF4}.

At every point in $\Omega_i$, moving frames are expressed in the following matrix form: 
\begin{align}
 \widehat{\mathbf{e}} =  {A} \widehat{\mathbf{x}},   \label{mf1}
\end{align}
where we introduced a new tensor $\widehat{\mathbf{e}} = \left [ \mathbf{e}^1,~\mathbf{e}^2,~ \mathbf{e}^3  \right ]^T$, $\widehat{\mathbf{x}} = \left [\mathbf{x}_1,~\mathbf{x}_2, ~\mathbf{x}_3 \right ]^T $ for the Cartesian coordinate unit vector $\mathbf{x}_i$, $1 \le i \le 3$. The matrix ${A}$ is known as the \textit{attitude matrix} \cite{ONeil} and represents the orientation of moving frames. By applying the differential operator for both sides, the 1-form $d \widehat{\mathbf{e} }$ is obtained as follows.
\begin{equation*}
d \widehat{\mathbf{e} } = \mathcal{W} \widehat{\mathbf{e}} ,
\end{equation*}
where $\mathcal{W}$ is a new tensor matrix representing $(dA) A^T $. The 1-form matrix $\mathcal{W}$, referred to as the \textit{connection form} \cite{Cartan1, Cartan2}, has nine components $[w^i_j]$ for $1 \le i,j \le 3$. Because of the orthonormality of moving frames, $\mathcal{W}$ is skew-symmetric and contains only three independent components as follows. 
\begin{equation*}
\mathcal{W}= \left [
\begin{array}{ccc}
0 & \omega^2_1 & \omega^3_1 \\
- \omega^2_1 & 0 &  \omega^3_2 \\
- \omega^3_1 & -\omega^3_2 & 0 
\end{array},
\right ] .
\end{equation*}
Because the component $\omega^i_j$ is also an 1-form, its value can be obtained when a specific direction is chosen. The connection form $\mathcal{W}$ can be obtained for the $\mathbf{e}^k$ direction, such as $\mathcal{W} \langle \mathbf{e}^k \rangle = d A \langle \mathbf{e}^k \rangle A^T$. Therefore,  $\omega^i_j \langle \mathbf{e}^k \rangle$ is obtained as follows.
\begin{equation} 
 \omega^i_j \langle \mathbf{e}^k \rangle =  ( {\mathbf{e}}^i )^T \cdot \mathbf{J}^j \cdot {\mathbf{e}}^k = \sum_{m=1}^3 e^i_{x_m} ( \nabla e^j_{x_m} \cdot \mathbf{e}^k  ) .  \label{omegaijk}
\end{equation}
For example, in the two-dimensional plane, various connection of moving frames can create a non-zero $\omega^2_1$ depending on the distribution of $\mathbf{e}^1$ for a certain direction, but $\omega^3_1$ and $\omega^3_2$ are zero regardless of the distribution of the moving frames. For a moving frame with unit length, the Christoffel symbol has the following relationship. $\mathbf{\mathbf{e}}^{\nu}_{;\alpha} =   \Gamma^{\mu}_{\alpha \nu}  \mathbf{e}^{\mu}$. Thus, this relationship reveals that $ \omega^i_j \langle \mathbf{e}^k \rangle$ is equivalent to $\Gamma^{\alpha}_{\mu \nu}$, i.e.,
\begin{equation} 
\Gamma^{\alpha}_{\mu \nu} = \omega^{\alpha}_{\mu} \langle \mathbf{e}^{\nu} \rangle   .  \label{ChristoffelConnection}
\end{equation}
Substituting Eq. \eqref{ChristoffelConnection} into Eq. \eqref{Covariant1}, we obtain 
\begin{equation}
\mathbf{u}_{;\alpha} = \left[ \nabla u^{\mu} \cdot \mathbf{e}^{\alpha} +  \omega^{\alpha}_{\mu} \langle \mathbf{e}^{\nu} \rangle u^{\nu} \right ] \mathbf{e}^{\mu}.  \label{CovariantMF}
\end{equation}
In comparison to Eq. \eqref{Covariant1}, the covariant formulation of Eq. \eqref{CovariantMF} neither requires the construction of curved axes nor the computation of $\Gamma^{\alpha}_{\mu \nu}$. Instead, moving frames are constructed at every point regardless of the underlying curvature of the domain. Moving frames are used as the direction derivative for scalar differentiation (first component) and the corresponding covariant compensation due to the changes in the axis (second component).

The first component of Eq. \eqref{Covariant1} is equivalent to the first component of Eq. \eqref{CovariantMF}. However, this is not true for the second component because the Christoffel symbol derived from Eq. \eqref{ChristoffelConnection} is derived from the axis with the unit tangent vector. If the Christoffel symbol is zero, or $g_{\alpha \beta}$ is constant, then Eq. \eqref{Covariant1} is equivalent to Eq. \eqref{CovariantMF}. For example, in the spherical coordinate axis on the sphere, Eq. \eqref{CovariantMF} yields only a low-order approximation to Eq. \eqref{Covariant1}. However, a \textit{special} construction of moving frames on a curved element can approximate Eq. \eqref{Covariant1} by Eq. \eqref{CovariantMF} with sufficiently high-order accuracy for the covariant derivative in the moving frames to function as a high-order method.

\begin{figure}[ht]
\centering
 \includegraphics[width=8cm]{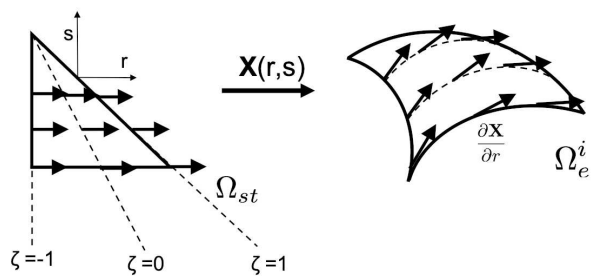} 
\caption{Mapping from the standard triangular element ($\Omega_{st}$) to a curved triangular element ($\Omega_{e}^i$).}
\label {Map}
\end{figure}

\section{LOCAL moving frames}

The computation of the covariant derivative in moving frames is exact if all metric tensors and Christoffel symbols are zero. A coordinate with this type of property is referred to as a \textit{Fermi coordinate} \cite{Misner}, and the corresponding moving frames are referred to as \textit{Euclidean}. However, it is nearly impossible to construct such a coordinate system on a generally curved surface, even on a sphere. In this paper, we introduce a convenient and efficient method for constructing moving frames to significantly reduce the geometric error caused by nontrivial metric tensors and Christoffel symbols.

In the finite element context, it is common to use a standard element ($\Omega_{st}$) for the mapping of a curved element, as depicted in Fig. \ref{Map}. A similar argument can be applied to the quadrilateral element, but we only focus on the mapping of a triangular element. Let $s$ and $r$ be the two Euclidean axes of the standard element in the range of $0 \le s,~r \le 1$. Let $\zeta$ be the another axis originating from one vertex, defined as $\zeta = 2(1+r)/(1-s)-1$ \cite{Spencerbook}. Let $\mathbf{X}(r,s)$ be the three-dimensional coordinate representation of the $i$th curved element $\Omega^i_e (x,y,z)$. Then, the differentiation of $\mathbf{X}(r,s)$ with respect to $\zeta$ produces the tangent vector $d \mathbf{X} / d \zeta$, which is nearly in the same direction as the tangent vector of the longitudinal axis $\theta$. $\Gamma^{\theta}_{\phi \phi}$ is zero for moving frames of unit length. If the first moving frame is aligned along the $\theta$ axis, then the error of $\Gamma^{\theta}_{\phi \phi}$ is equivalent to $-2 \sin 2 \theta$, implying that the the error of the covariant derivative is the first order of convergence with respect to the length of the edge $\ell$, i.e., $\mathcal{O} (\ell)$.

The other option is to construct the moving frames parallel to each edge, or along the axis of $s$ and $r$, respectively. The differentiation of $\mathbf{X}(r,s)$ with respect to the $s$-axis produces the tangent vector $d \mathbf{X} / d s$ parallel to the line of $r$ = $\mbox{constant}$. A similar argument can be applied to $d \mathbf{X} / d r$ but, in general, $(d \mathbf{X} / d r) \cdot (d \mathbf{X} / d s) \neq 0$. Contrary to $d \mathbf{X} / d \zeta$, $d \mathbf{X} / d s$ is almost Euclidean in the element, and its orthonormal vector is also almost Euclidean. Even though the constructed frames are not the exact Fermi coordinate system, they are sufficiently Euclidean for significantly reduced corresponding error. Let us refer these frames as \textit{LOCAL moving frames}, whereas the moving frames aligned along the spherical coordinate axis are referred to as the \textit{spherical} moving frames. 

LOCAL moving frames can be easily constructed as follows. Consider the three edges ($1 \le E_k \le 3$) of a curved triangular. Moving frames are constructed along each edge such that $e_{E_1} = d \mathbf{X} / d r$, $e_{E_2} = d \mathbf{X} / d s$, and $e_{E_3} = 0.5 ( d \mathbf{X} / d r  + d \mathbf{X} / d s )$. Then, the LOCAL moving frames with the lowest magnitude of the covariant divergence are chosen as follows. 
\begin{equation*}
\mathbf{e} = \mathbf{e}_{E_k}~ \mbox{for the index}~E_k ~\mbox{corresponding to} ~\min_{k=1}^3 \left \{ \left \| \sum_i^2 \nabla \cdot \mathbf{e}_{E_k}^i \right \| \right \}, 
\end{equation*}

LOCAL moving frames are similar to spherical moving frames around the equator because moving frames are mostly Euclidean in those regions. However, LOCAL moving frames are generally discontinuous across the elements, contrary to spherical moving frames.

\begin{figure}[ht]
\centering
\subfloat[$\nabla \cdot \mathbf{e}^1_{sph}$ ] {\label{sphe1}  \includegraphics[width=3cm]{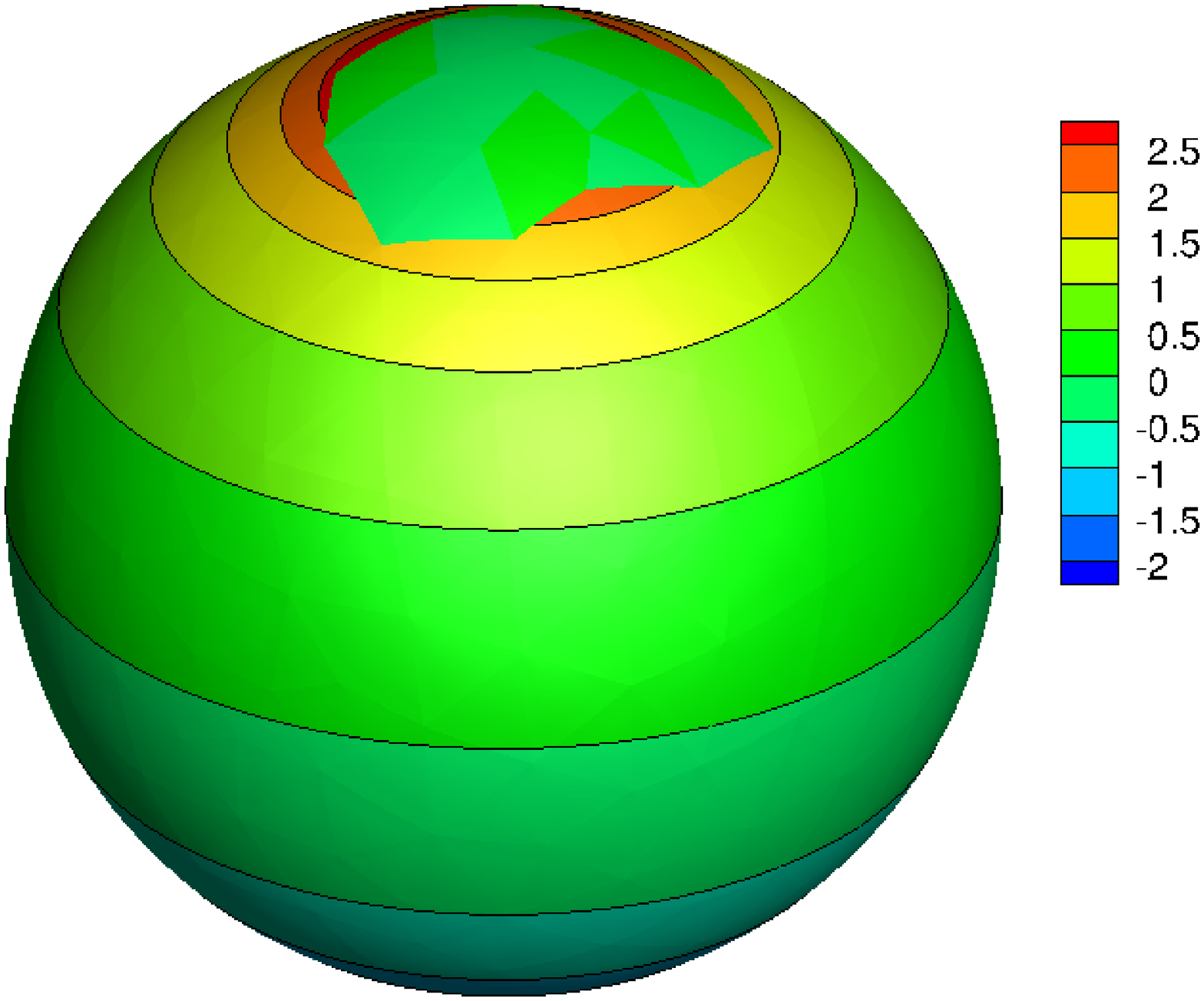}  }
\subfloat[$\nabla \cdot \mathbf{e}^2_{sph}$ ] {\label{sphe2}  \includegraphics[width=3cm]{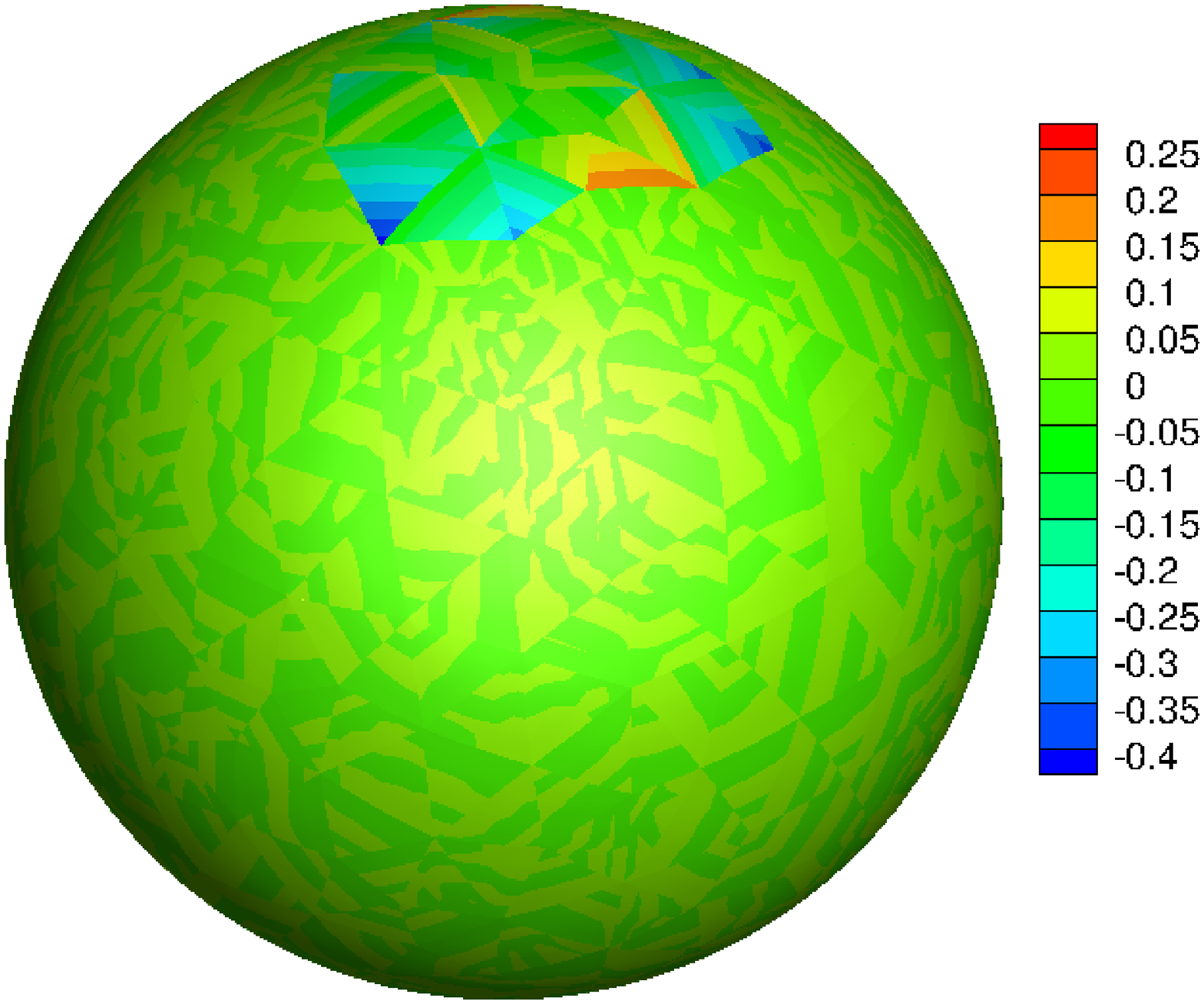}  }
\subfloat[$\nabla \cdot \mathbf{e}^1_{loc}$ ] {\label{loce1}  \includegraphics[width=3cm]{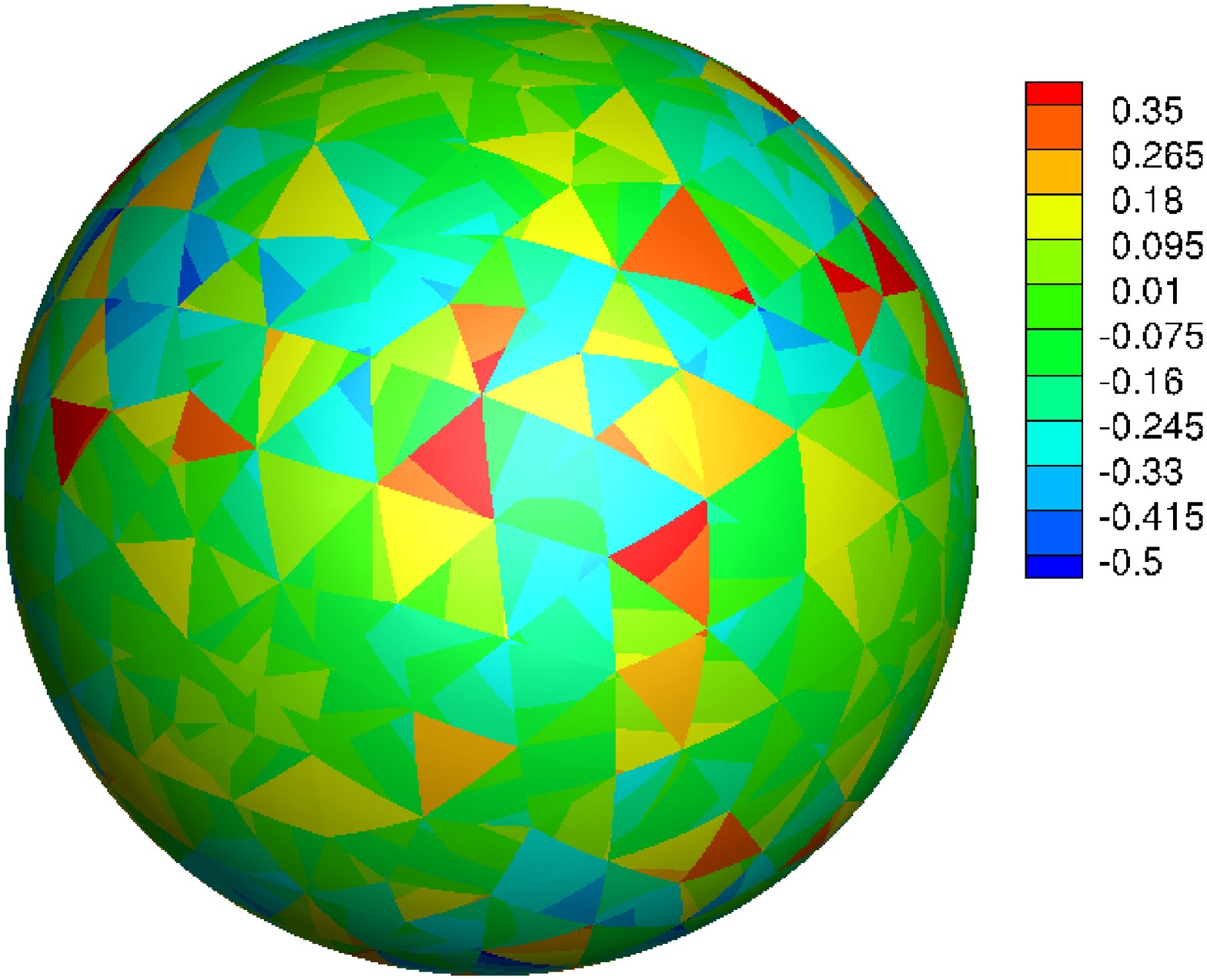}  }
\subfloat[$\nabla \cdot \mathbf{e}^2_{loc}$ ] {\label{loce2}  \includegraphics[width=3cm]{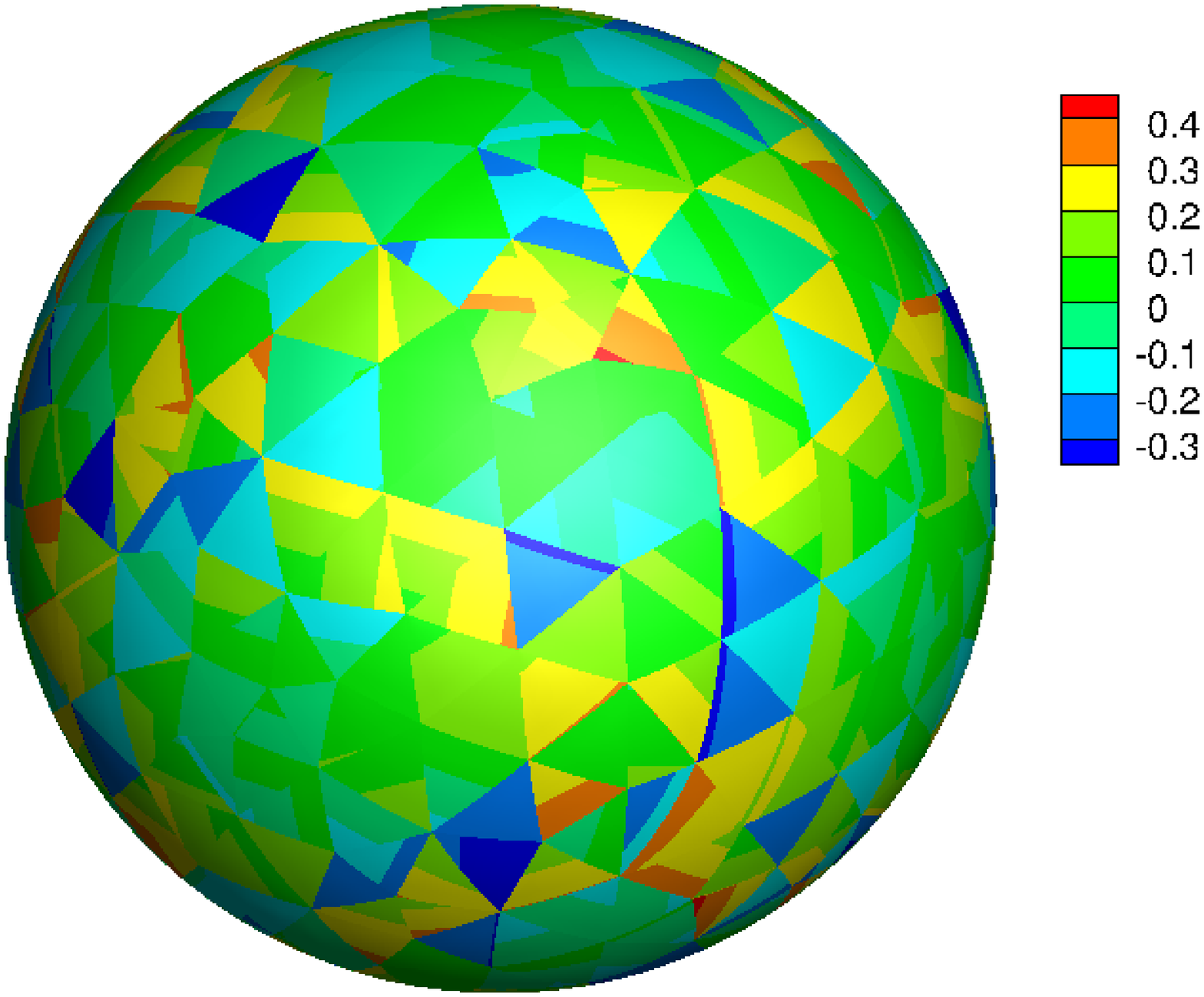}  }
\caption{ Distribution of divergence of moving frames $\mathbf{e}^i$ and $\mathbf{e}^i$ for spherical ($\mathbf{e}^i_{sph}$) and LOCAL moving frames ($\mathbf{e}^i_{loc}$) .}
\label {LOCALmf}
\end{figure}

Fig. \ref{LOCALmf} presents the difference of $\nabla \cdot \mathbf{e}^i,~ i=1,2$ between spherical moving frames and LOCAL moving frames. Consider a tessellated spherical mesh with 498 elements and a $4.99202e$-$8$ mesh error. For spherical moving frames, $\nabla \cdot \mathbf{e}^1$ and $\nabla \cdot \mathbf{e}^2$ are $0.774803$ and $0.0425204$, respectively, whereas for LOCAL moving frames, $\nabla \cdot \mathbf{e}^1$ and $\nabla \cdot \mathbf{e}^2$ are $0.192841$ and $0.16384$, respectively. The strategy of this scheme is to lower the maximum of $\nabla \cdot \mathbf{e}^i$ and distribute it equally to both moving frames, $\mathbf{e}^1$ and $\mathbf{e}^2$. This reconstruction of moving frames reduces the geometric error significantly compared to the reduction of the discretization error in differentiation and integration. In the next section, we will demonstrate that the derived moving frames with the connection form significantly increases the accuracy of the covariant derivative.

\section{Covariant formulation and test cases}

Consider a unit sphere with the following metric $ds^2 = r^2 d \theta^2 + r^2 \sin^2 \theta d \phi^2 $. The velocity field of the Rossby-Haurwitz wave, popular in the shallow water equations, is defined as
\begin{align*}
\mathbf{v} &= v_{\phi}  \hat{\boldsymbol{\phi}} + v_{\theta}  \hat{\boldsymbol{\theta}} ,\\
\mbox{where}&, \\
v_{\phi} &=  \omega \sin \theta + K \sin^3 \theta ( 4 \cos^2 \theta - \sin^2 \theta ) \cos 4 \phi, \\
v_{\theta} &=  -4 K  \sin^3 \theta \cos \theta \sin 4 \phi , 
\end{align*}
where $\omega = K = 7.848 \times 10^{-6}~ s^{-1}$. In LOCAL moving frames, the vector $\mathbf{v}$ is expanded as $\mathbf{v} = v_1 \mathbf{e}^1 + v_2 \mathbf{e}^2$ for almost Euclidean moving frames $\mathbf{e}^1$ and $\mathbf{e}^2$.

\subsection{Gradient}
The first test relies on the fact that the gradients of a scalar variable on a curved surface should be equal, independent of the axis. For a spherical axis of $(\theta, \phi)$, the gradient of a scalar variable, such as, $v_{\phi}$, is given as
\begin{equation}
\nabla v_{\phi}  = \frac{\partial v_{\phi} }{\partial \theta} \boldsymbol{\theta} + \frac{1}{\sin^2 \theta} \frac{\partial v_{\phi} }{\partial \phi} \boldsymbol{\phi} , \label{gradf}
\end{equation}
where
\begin{align*}
\frac{\partial v_{\phi} }{\partial \theta} &= \omega \cos \theta + K \sin^2 \theta [ 3 \cos \theta ( 4 \cos^2 \theta - \sin^2 \theta ) - 10 \sin^2 \theta \cos \theta ] \cos 4 \phi, \\
\frac{\partial v_{\phi} }{\partial \phi} &= - 4 K \sin \theta ( 4 \cos^2 \theta - \sin^2 \theta ) \sin 4 \phi .
\end{align*}
The computation of the gradient in LOCAL moving frames should have the same value, expressed as follows,
\begin{equation}
\nabla v_{\phi}  = ( \nabla v_{\phi}  \cdot \mathbf{e}^1 ) \mathbf{e}^1 +  ( \nabla v_{\phi}  \cdot \mathbf{e}^2 ) \mathbf{e}^2   .  \label{CovGrad}
\end{equation}
Computationally, this implies that Eq. \eqref{gradf} should converges to Eq. \eqref{CovGrad} as $p$ increases or $h$ decreases. Table \ref{Gradtest} confirms the the exponential convergence of the difference between the two formulations for the gradient of $v_{\phi}$.

\begin{table}
\caption{Difference between Eq. \eqref{gradf} and Eq. \eqref{CovGrad}. Sphere of radius $1.0$. $h=0.4$ and 498 elements. Several elements close to the poles are not considered due to the singularities of the spherical coordinate axis.}
\begin{center}
\begin{tabular}{|c|c|c|c|c|c|c|}
\hline 
   p  &   3     &      4       &     5   &     6  &  7   &  8 \\
   \hline
  Diff & 1.85e-4 & 8.64e-06 & 3.75e-07 & 2.21e-08 & 1.02e-09 & 6.62e-11 \\
 \hline
\end{tabular} 
\end{center}
\label{Gradtest}
\end{table}%

\subsection{Divergence}
In the spherical coordinate axis, the divergence of the velocity vector $\mathbf{v}$ is obtained by the following covariant formulation.
\begin{equation}
\nabla \cdot \mathbf{v} = \frac{1}{\sin \theta} \left ( \frac{\partial v_{\phi}}{\partial \phi} + \frac{\partial }{\partial \theta} ( v_{\theta} \sin \theta )  \right ) \label{DivSph},
\end{equation}
In moving frames, the divergence is obtained as follows.
\begin{equation}
\nabla \cdot \mathbf{v} =  \nabla v_1 \cdot  \mathbf{e}^1 + \Gamma^1_{21} v_2 + \nabla v_2 \cdot  \mathbf{e}^2 + \Gamma^2_{21} v_1,
 \label{DivDirect}
\end{equation}
where the Christoffel symbol $\Gamma^i_{jk}$ is computed by the connection 1-form of $\omega^i_j \langle \mathbf{e}^k \rangle$, as depicted in Eq. \eqref{ChristoffelConnection}. The divergence of the velocity vector of the Rossby-Haurwitz wave is analytically zero.

\begin{figure}[ht]
\centering
\subfloat[Divergence] {\label{Divpconv}  \includegraphics[width=5cm]{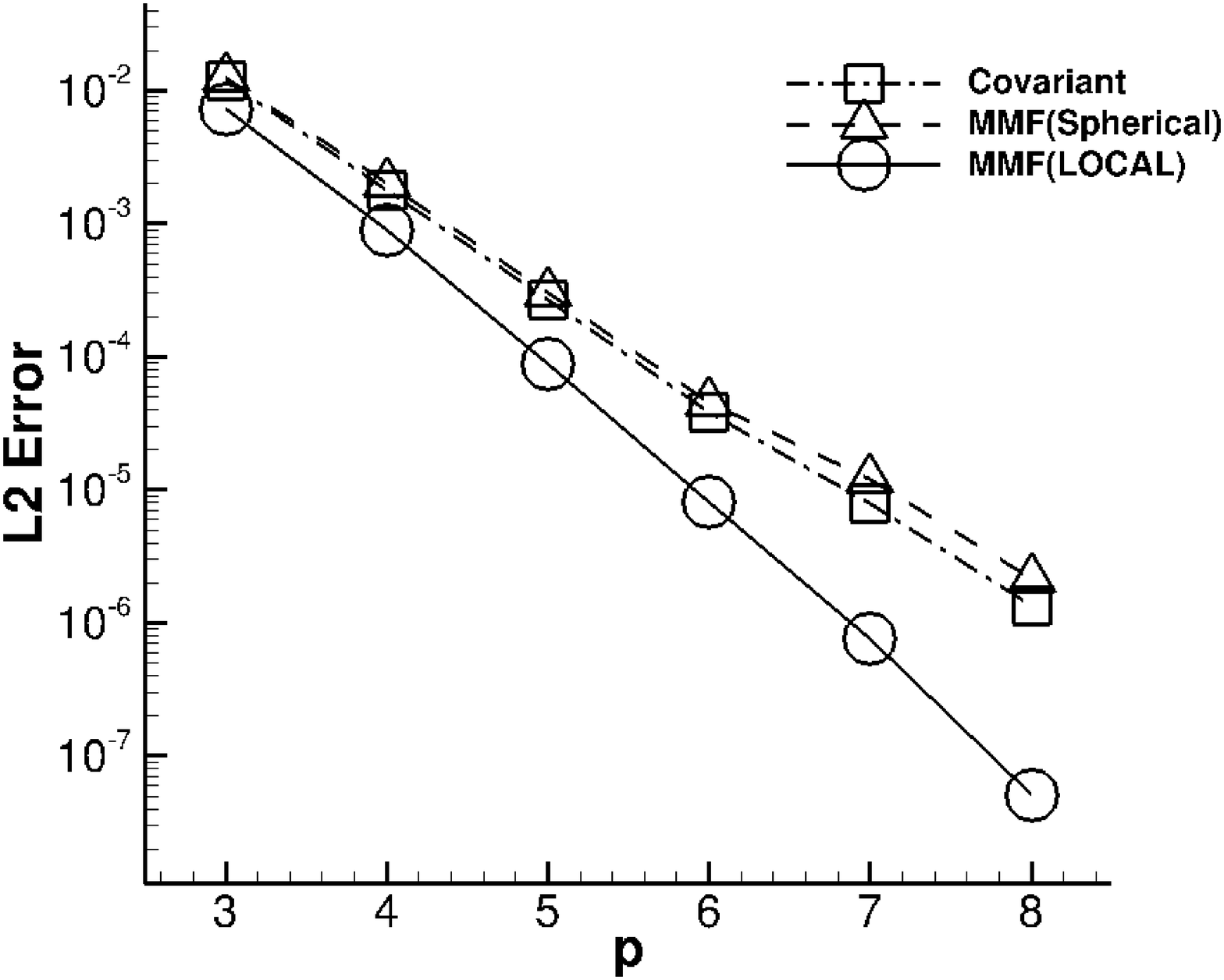}  }
\subfloat[Curl] {\label{Curlpconv}  \includegraphics[width=5cm]{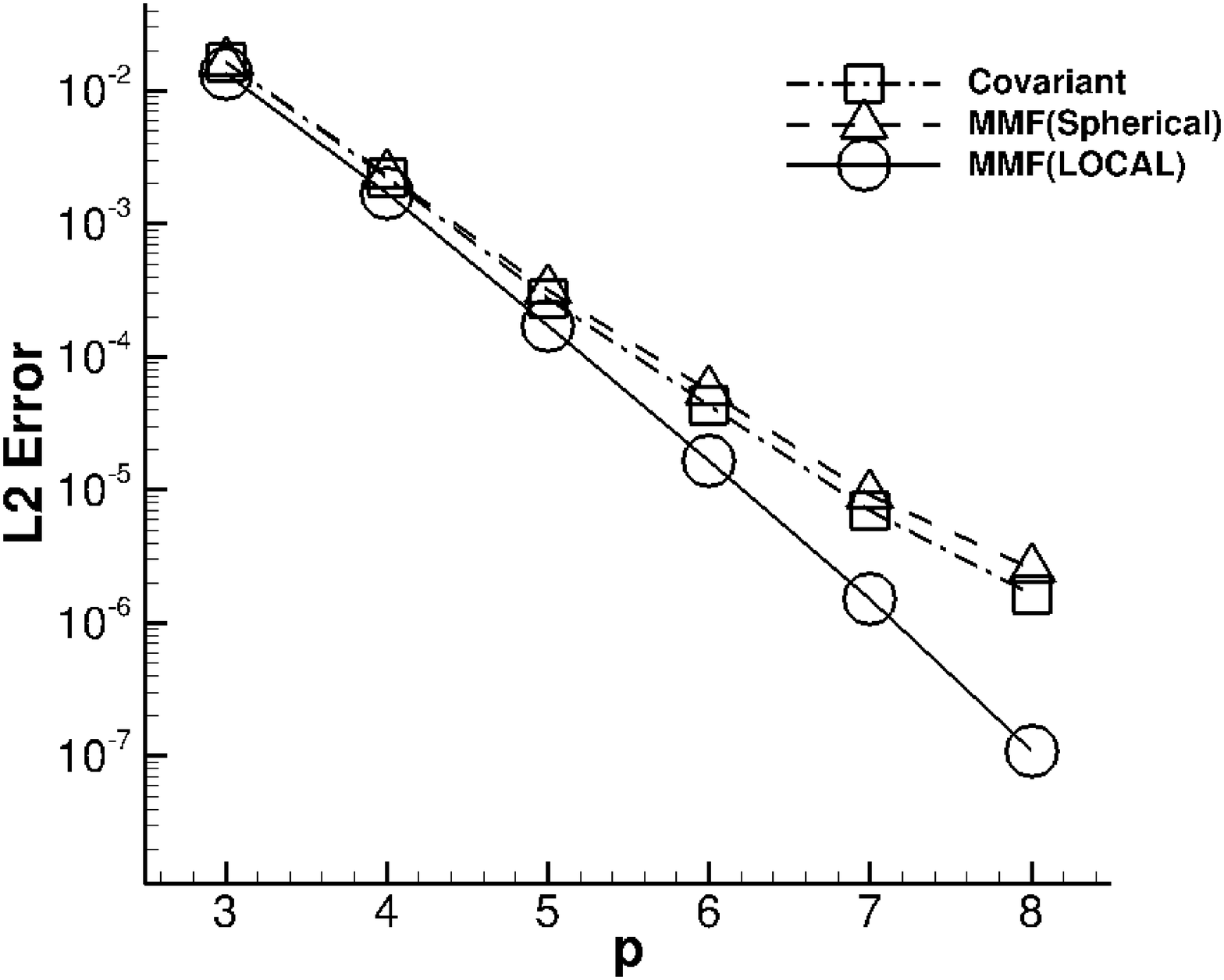}  }
\caption{$p$-convergence of divergence and curl on the sphere. $h$=$0.3445$.}
\label {pconv}
\end{figure}

Fig. \ref{Divpconv} illustrates the exponential convergence by Eq. \eqref{DivSph} (Covariant), Eq. \eqref{DivDirect} with spherical moving frames (MMF (Spherical)), and Eq. \ref{DivDirect} with LOCAL moving frames (MMF(LOCAL)). Fig. \ref{Divpconv} confirms that Eq. \eqref{DivDirect} with LOCAL moving frames has the highest accuracy. The difference between the other method becomes larger as $p$ increase, which implies that, as $p$ increases, the geometric error contributes more to the overall error. Table \ref{Divtesthconv} presents the convergence order for the three methods to demonstrate that Eq. \eqref{DivDirect} with LOCAL moving frames is the most accurate with an improved convergence order. The order should be theoretically $p$ for $p$=$5$ because it is the first derivative of a vector. However, the geometric error of the mesh undermines this order, ending up with $3.47,~ 5.3,~ 4.03$ for covariant formulation and $3.27,~ 5.25,~ 4.04$ for spherical moving frames. For LOCAL moving frames, the order is increased to be nearly equivalent to the ideal order: $5.08,~4.79,~5.58$.

\begin{table}
\caption{$h$-convergence of divergence on the sphere by the three different methods. $p=5$. }
\begin{center}
\begin{tabular}{|c|c|c|c|c|}
\hline 
   h  &   0.186726       &        0.262293        &      0.3445      &      0.502745 \\
   $N_e$  & 1918    &    970   &  498     &      278 \\
\hline \hline 
   Covariant  & 1.93009e-05       &      6.28588e-05     &     0.000268214      &      0.0012311  \\
 order & -  & 3.4746  &  5.3218  &  4.0315  \\
    MMF (Sphere)  & 2.38394e-05      &        7.23155e-05       &     0.000302484      &     0.00139154  \\
 order & -  &  3.2655  &  5.2488  &  4.0375 \\
     MMF (LOCAL)  & 4.26305e-06   &   2.3937e-05   &    8.82688e-05    &  0.000728567  \\
 order & -  &  5.0775  & 4.7865   & 5.5840 \\
 \hline
\end{tabular} 
\end{center}
\label{Divtesthconv}
\end{table}%

\subsection{Curl}
For the computation of $\mathbf{k} \cdot \nabla \times \mathbf{v}$ for the surface normal vector $\mathbf{k}$, the covariant computation of the curl operator in the spherical coordinate axis is given as
\begin{equation}
\mathbf{k} \cdot ( \nabla \times \mathbf{v} ) =  \frac{1}{\sin \theta} \left ( \frac{\partial v_{\theta}}{\partial \phi} - \frac{\partial }{\partial \theta} ( v_{\phi} \sin \theta )  \right ) \label{CurlSph},
\end{equation}
By direct differentiation in moving frames, the curl can be computed as
\begin{equation}
\mathbf{k} \cdot (  \nabla \times  \mathbf{v} ) =  \nabla v_2 \cdot  \mathbf{e}^1 + \Gamma^2_{12} v_2 - (  \nabla v_1 \cdot  \mathbf{e}^2 + \Gamma^1_{21} v_1 ).
 \label{CurlMF}
\end{equation}
The analytical value of $\mathbf{k} \cdot (  \nabla \times  \mathbf{v} )$ for the velocity vector of the Rossby-Haurwitz wave is given as
\begin{equation*}
(\nabla \times \mathbf{v} ) \cdot \mathbf{r} = - 2 \omega \cos \theta + 30 K \sin^4 \theta \cos \theta \cos 4 \phi .
\end{equation*}
Fig. \ref{Curlpconv} illustrates the similar exponential convergence as that of the divergence: Eq. \eqref{CurlSph} (Covariant), Eq. \eqref{CurlMF} with spherical moving frames (MMF (Spherical)), and Eq. \ref{CurlMF} with LOCAL moving frames (MMF(LOCAL)). Similarly, Eq. \eqref{CurlMF} with LOCAL moving frames has the highest accuracy. Table \ref{Curltesthconv} presents the convergence order for the three methods, which indicates that Eq. \eqref{CurlMF} with LOCAL moving frames exhibits an improved convergence order of $5.0775,~ 4.7865 ,~ 5.5840$, closer to the ideal spectral convergence of $p$.

\begin{table}
\caption{$h$-convergence of curl on the sphere by the three methods. $p=5$. }
\begin{center}
\begin{tabular}{|c|c|c|c|c|}
\hline 
   h  &   0.186726       &        0.262293        &      0.3445      &      0.502745 \\
   $N_e$  & 1918    &    970   &  498     &      278 \\
\hline \hline 
   Covariant  & 1.95592e-05  &  6.83048e-05  &   0.000277284 &  0.00141376 \\
 order & -  & 3.4746  &  5.3218  &  4.0315  \\
    MMF (Sphere)  & 1.92204e-05   & 7.45893e-05   &   0.00031796   &  0.00171716  \\
 order & -  &  3.2655  &  5.2488  &  4.0375 \\
     MMF (LOCAL)  & 7.66923e-06  &  4.07124e-05   &  0.000173096    &.  0.00108435  \\
 order & -  &  5.0775 &   4.7865  &  5.5840 \\
 \hline
\end{tabular} 
\end{center}
\label{Curltesthconv}
\end{table}%

\section{Helmholtz-Hodge Decomposition}
On a curved surface $\Omega$ with Neumann boundary or no boundary, the Helmholtz-Hodge decomposition (HHD) finds the unique three components of a vector field $\mathbf{v}$, similar to \cite{Bhatia}
\begin{equation}
\boldsymbol{\xi} = \nabla u + \nabla \times \mathbf{v} + \mathbf{h},  \label {HHD}
\end{equation}
where $\nabla u$ is a curl-less irrotational vector, $\nabla \times \mathbf{v}$ is a divergence-less incompressible vector, and $\mathbf{h}$ is a harmonic vector with zero vector Laplacian, i.e., $\nabla^2 \mathbf{h} = \mathbf{0}$. The irrotational component is obtained by applying the divergence to Eq. \eqref{HHD}. The incompressible component is first expressed as $\nabla \times \mathbf{R} = J \nabla R$ for the linear operator $J$, transforming $\mathbf{v} = v^1 \mathbf{e}^1 + v^2 \mathbf{e}^2$ into $J \mathbf{v} = -v^2 \mathbf{e}^1 + v^1 \mathbf{e}^2$, and is obtained by applying the divergence to Eq. \eqref{HHD}, i.e.,
\begin{align}
\nabla^2 u &= \nabla \cdot \boldsymbol{\xi} - \frac{1}{A} \int \nabla \cdot \boldsymbol{\xi} dx,     \label{Helm1} \\ 
\nabla^2 v &= - \nabla \cdot J \boldsymbol{\xi} + \frac{1}{A} \int \nabla \cdot J \boldsymbol{\xi} dx, \label{Helm2}
\end{align}
where $A$ is the surface area of the domain $\Omega$. The second components in the right-hand side is added because the domain has a Neumann boundary condition or no boundaries. Then, the vector $\mathbf{h}$ is obtained by subtracting the two components from $\mathbf{v}$, i.e., $\nabla \cdot \mathbf{h} = ( \int \nabla \cdot \boldsymbol{\xi} dx ) / A$, $\nabla \cdot J \mathbf{h} = ( \int \nabla \cdot J \boldsymbol{\xi} dx ) / A$. Because these values are constant in the domain, the vector Laplacian of $\nabla^2 \mathbf{h} = \nabla (\nabla \cdot \mathbf{h} ) + \nabla \times (\nabla \times \mathbf{h})$ is zero.

On a surface, Eqs. \eqref{Helm1} and \eqref{Helm2} are covariant derivatives, which should be computed by Eq. \eqref{DivDirect} and Eq. \eqref{CurlMF}, respectively. Inaccurate computation on the right-hand side of Eq. \eqref{Helm1} and \eqref{Helm2} fail to locate the exact source of the flow represented as the irrotational and incompressible components. The scheme is implemented at the open-source spectral/hp library, referred to as \textit{Nektar++} \cite{Nektar++}. Eqs. \eqref{Helm1} and \eqref{Helm2} are solved by the built-in Helmholtz solver in the context of continuous or discontinuous Galerkin methods.

\begin{figure}[ht]
\centering
\subfloat[Irrotational] {\label{SphereIrr}  \includegraphics[width=4cm]{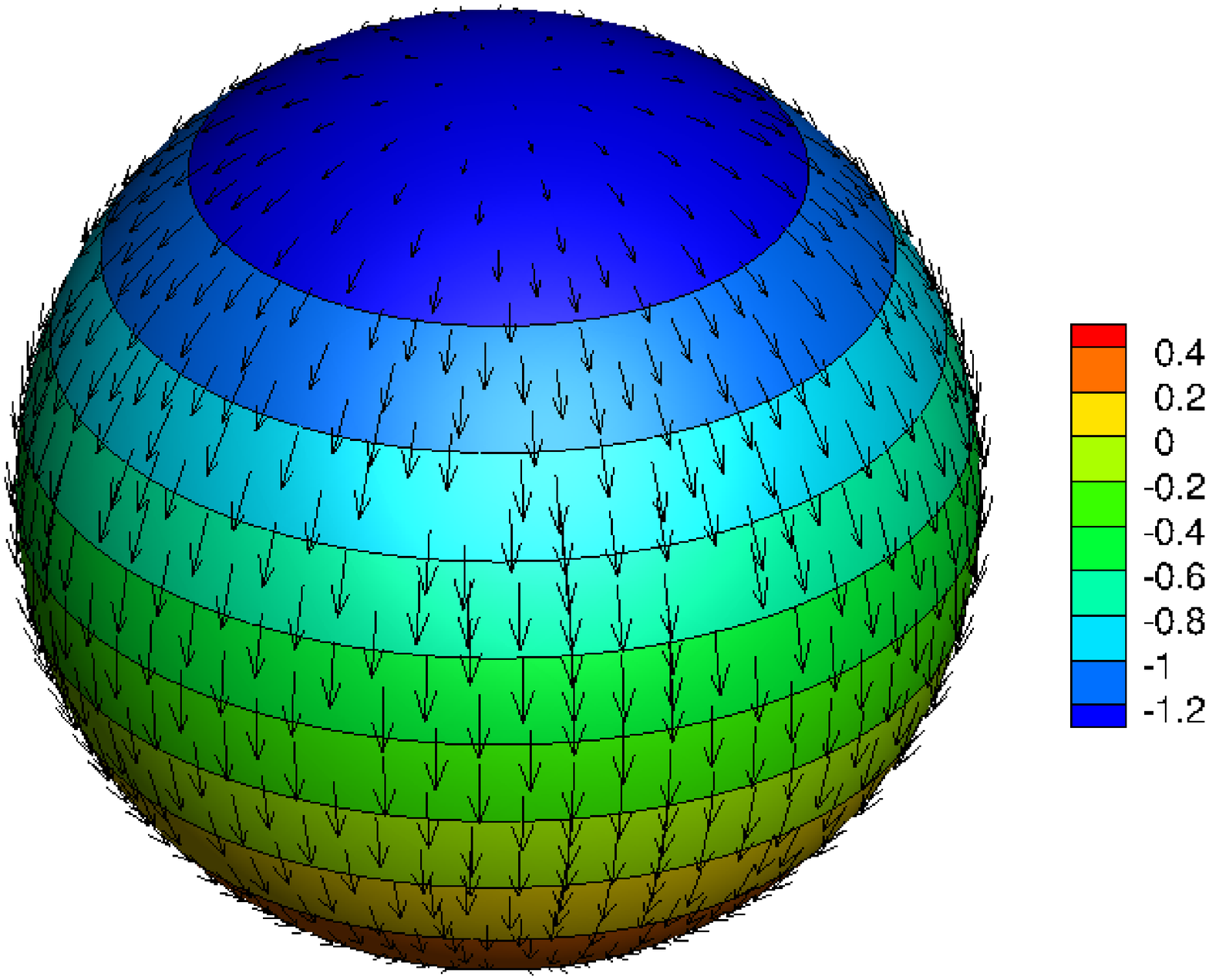}  }
\subfloat[Incompressible] {\label{SphereIncomp}  \includegraphics[width=4cm]{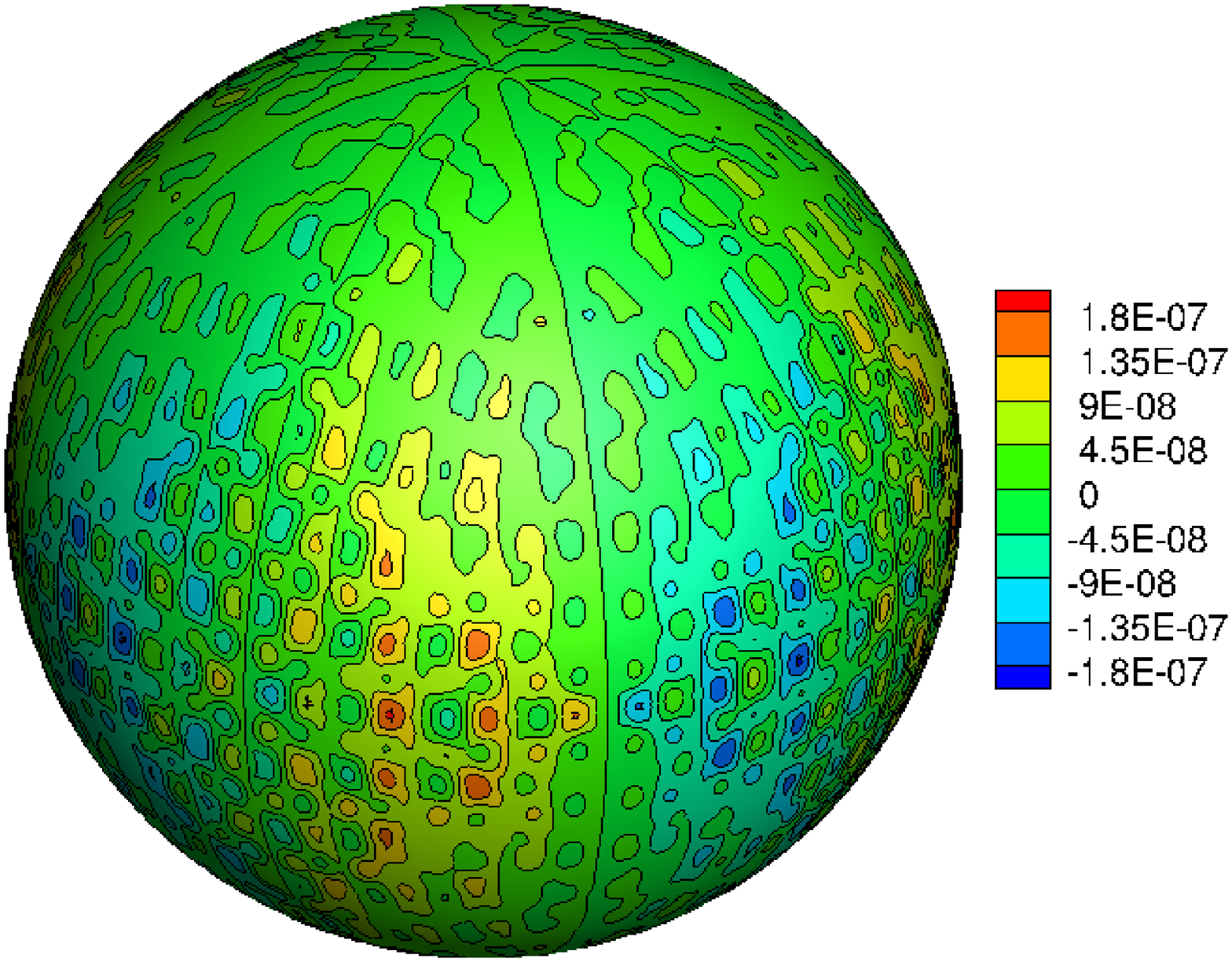}  }
\subfloat[Harmonic] {\label{Sphereharmon}  \includegraphics[width=4cm]{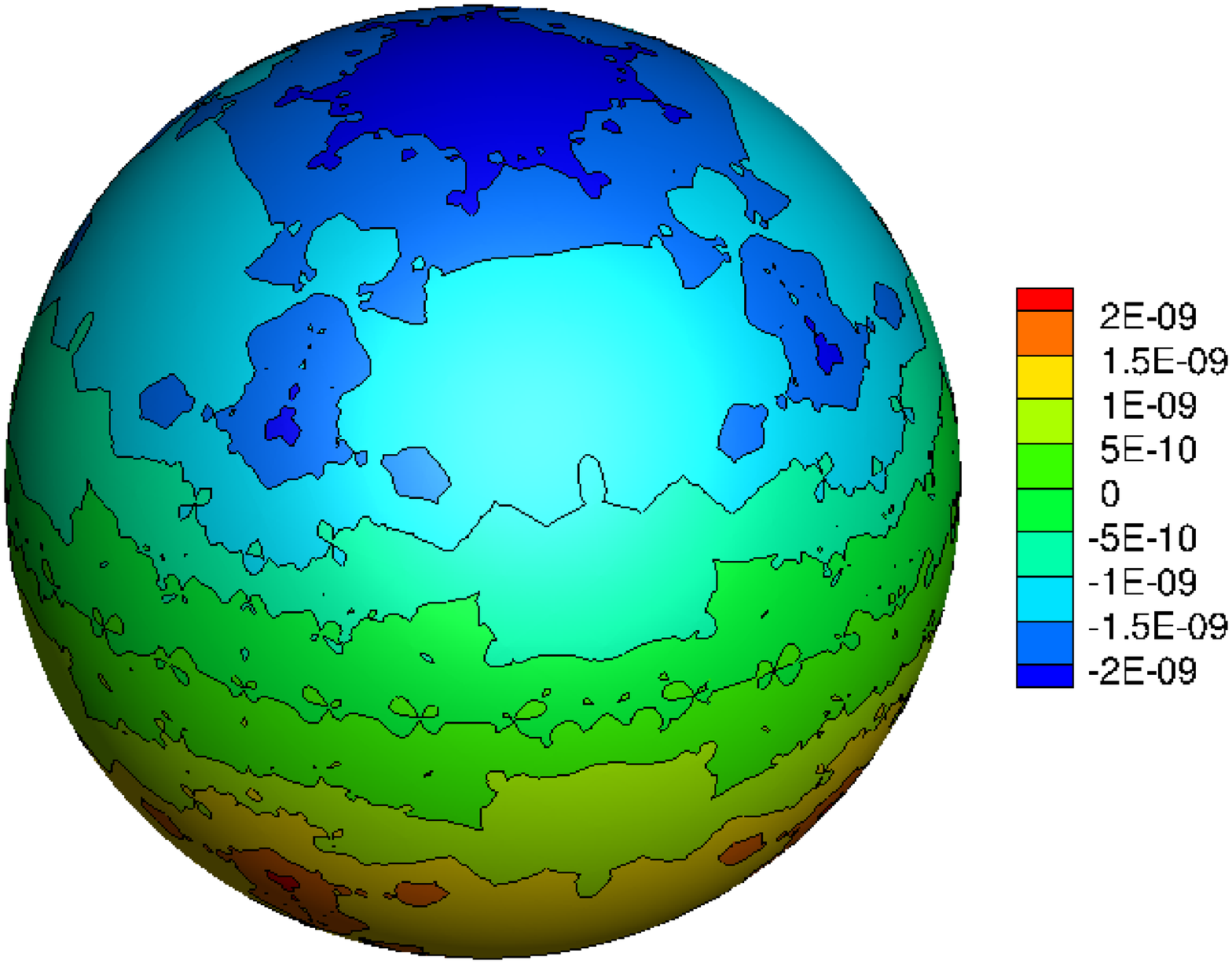}  }
\caption{HHD of a curl-less vector, the spherical moving frames multiplied by $\sin \theta$. (a) The irrotational component of the vector with the potential $u$, (b) the potential of the incompressible flow $v$, and (c) the potential of harmonic flow $U$ such as $\mathbf{h} = \nabla U$. $h$=$0.2$. $p$=$10$.}
\label {HHDsphere}
\end{figure}

\begin{figure}[ht]
\centering
\subfloat[Irrotational] {\label{RHIrr}  \includegraphics[width=4cm]{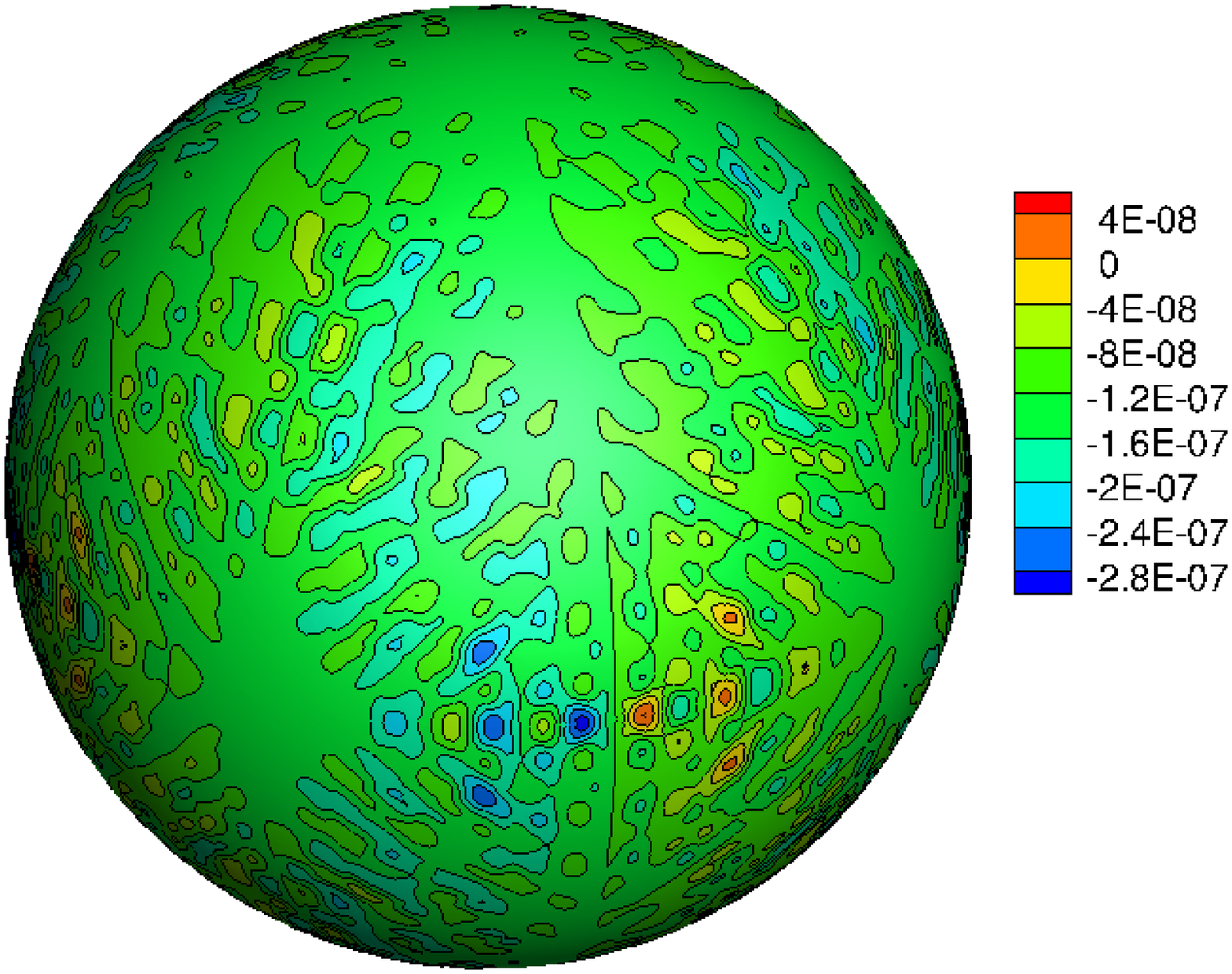}  }
\subfloat[Incompressible] {\label{RHIncomp}  \includegraphics[width=4cm]{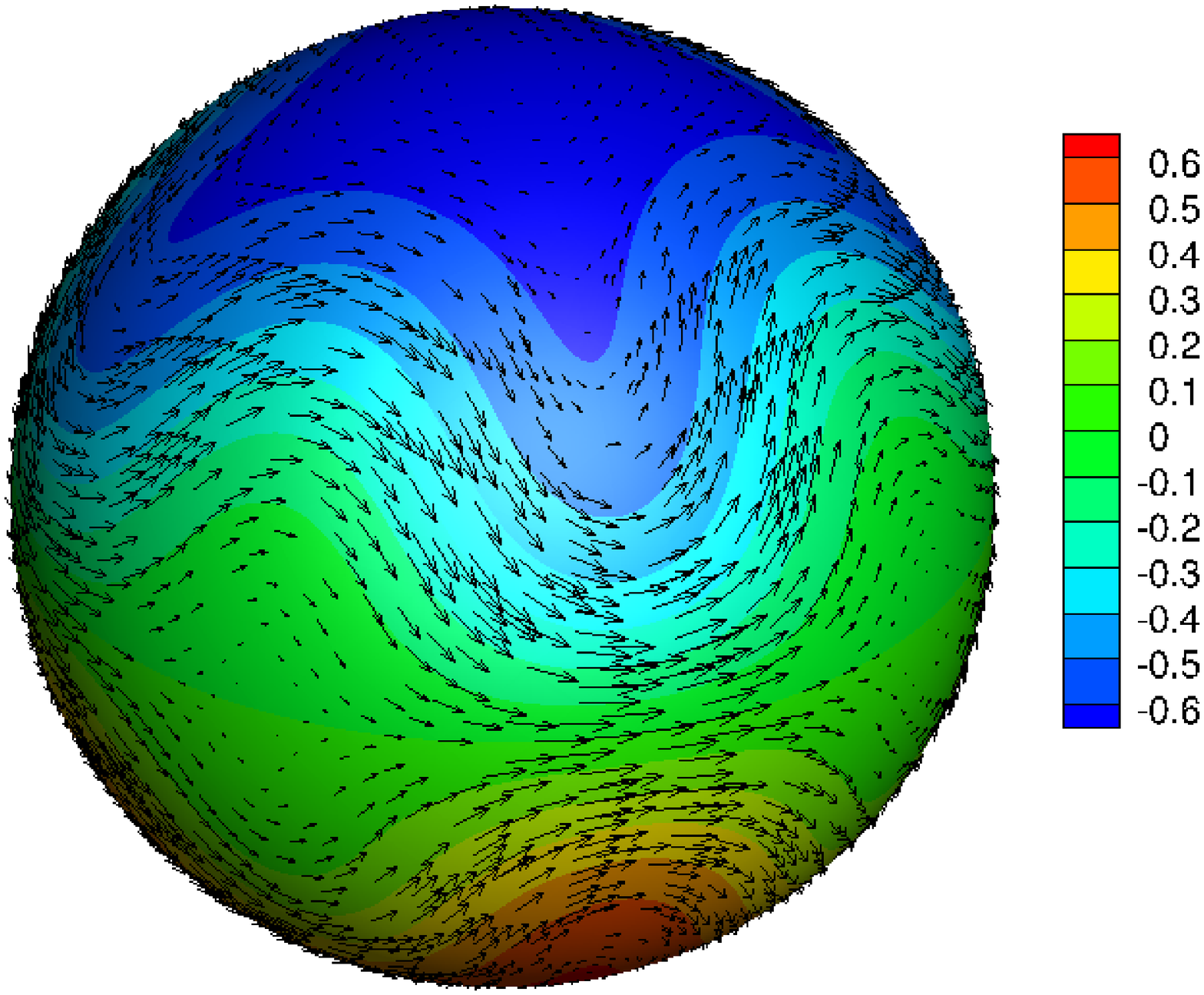}  }
\subfloat[Harmonic] {\label{RHharmon}  \includegraphics[width=4cm]{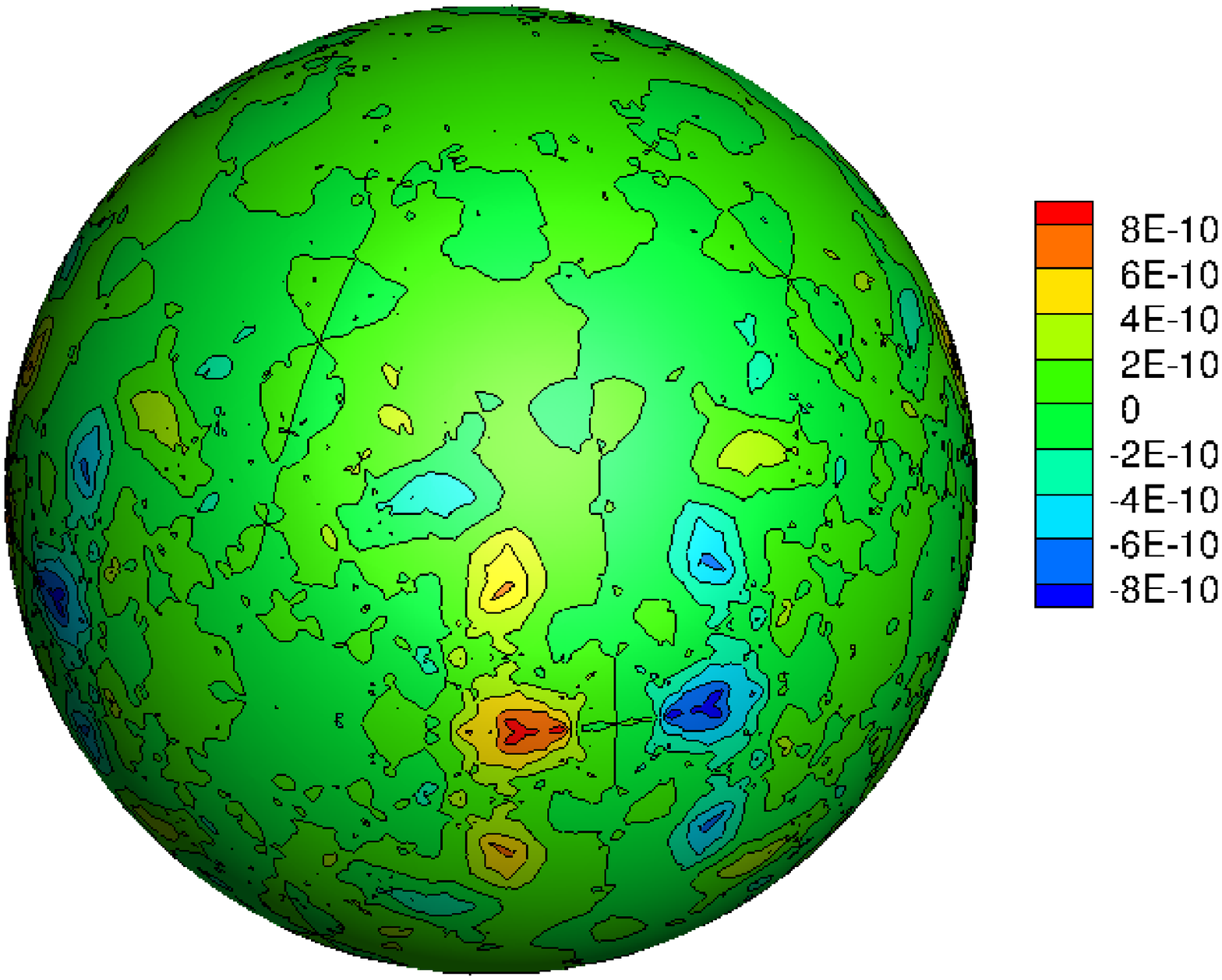}  }
\caption{HHD of a divergence-less vector, the Rossby-Haurwitz velocity vector. (a) The potential $u$ of the irrotational component, (b) the potential $v$ of the incompressible flow, and (c) the potential $U$ of harmonic flow such as $\mathbf{h} = \nabla U$. $h$=$0.2$. $p$=$10$.}
\label {HHDRHwave}
\end{figure}

Fig. \eqref{HHDsphere} represents the HHD of a curl-less vector with the following error: $\| \nabla \times \nabla u \| = 2.44e$-$9$, $\| \nabla \cdot  J \nabla v \| = 2.77e$-$9$. $\| \nabla \cdot \mathbf{h} -  \frac{1}{A} \int \nabla \cdot \boldsymbol{\xi} dx \| = 6.60e$-$6$, $\| \nabla \cdot J \mathbf{h} -  \frac{1}{A} \int \nabla \cdot J \boldsymbol{\xi} dx \| = 6.51e$-$6$, and $\| \nabla^2 \mathbf{h} \| = 1.79e$-$3$. Fig. \eqref{HHDRHwave} illustrates the HHD of the divergence-less Rossby-Haurwitz velocity vector with the following error. $\| \nabla \times \nabla u \| = 1.95e$-$9$, $\| \nabla \cdot  J \nabla v \| = 2.40e$-$9$. $\| \nabla \cdot \mathbf{h} -  \frac{1}{A} \int \nabla \cdot \boldsymbol{\xi} dx \| = 6.16e$-$6$, $\| \nabla \cdot J \mathbf{h} -  \frac{1}{A} \int \nabla \cdot J \boldsymbol{\xi} dx \| = 5.51e$-$6$, and $\| \nabla^2 \mathbf{h} \| = 1.63e$-$3$.

\begin{figure}[ht]
\centering
\subfloat[Irrotational] {\label{AtriumIrrot}  \includegraphics[width=4cm]{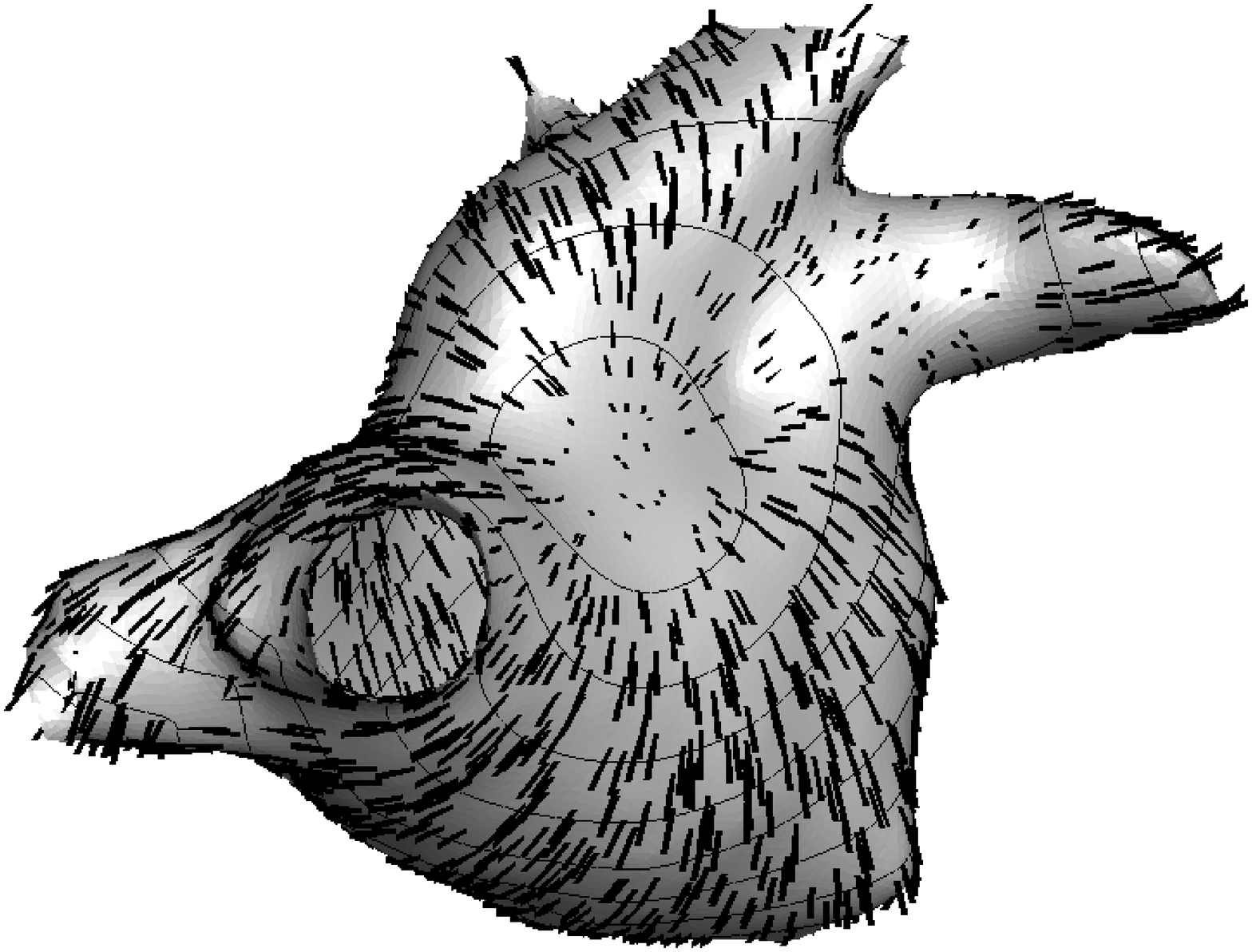}  }
\subfloat[Incompressible] {\label{Atrium}  \includegraphics[width=4cm]{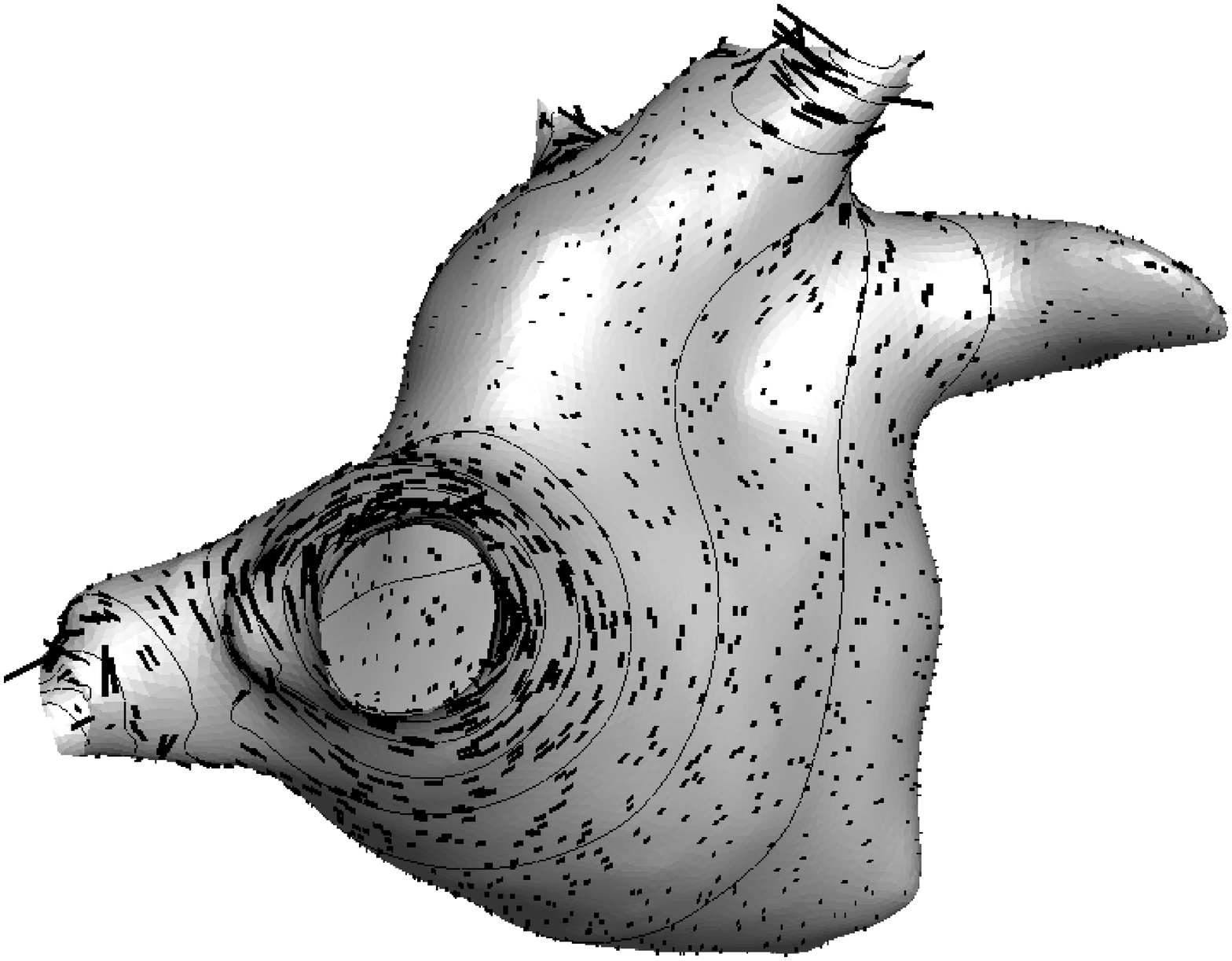}  }
\subfloat[Harmonic] {\label{AtriumIncomp}  \includegraphics[width=4cm]{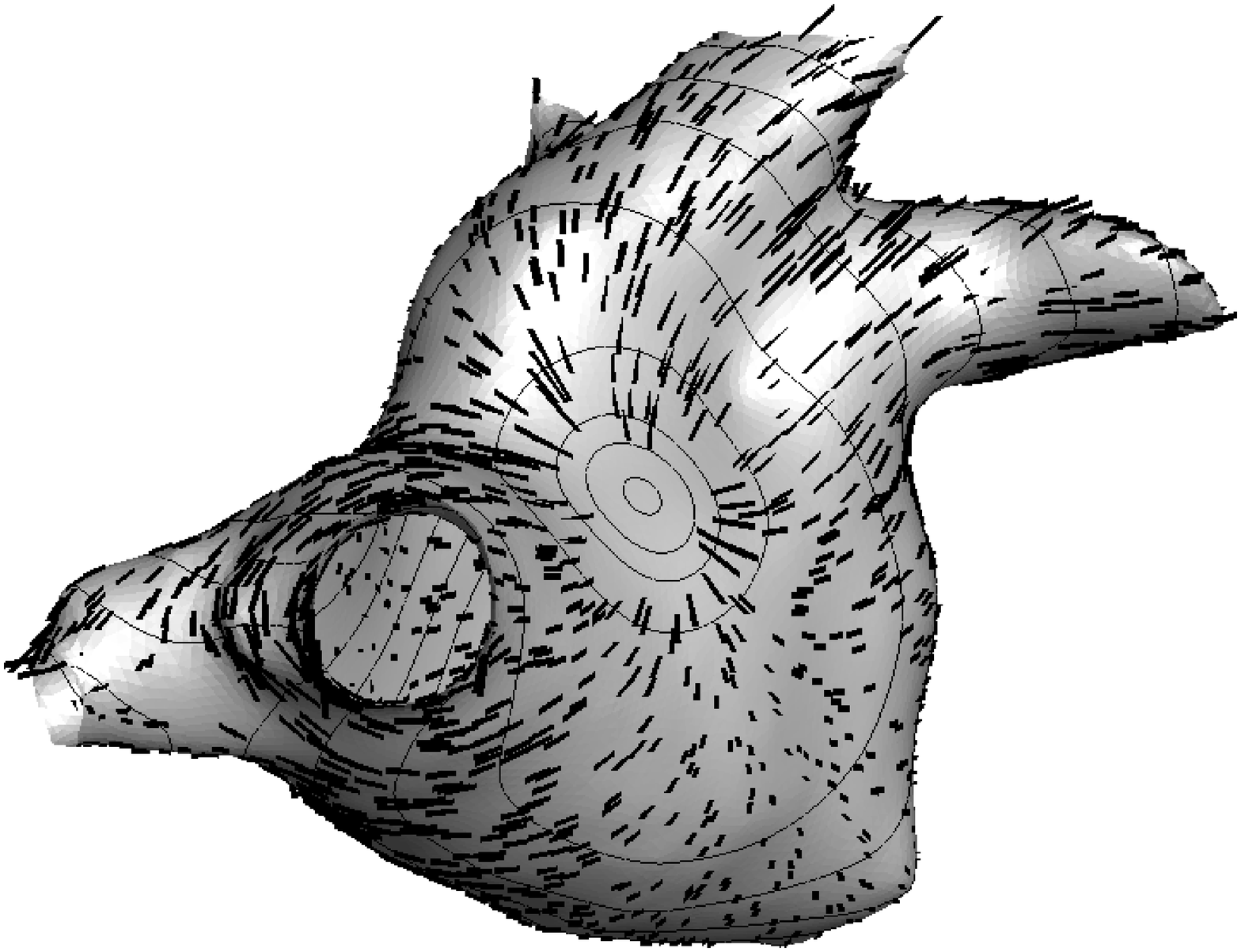}  }
\caption{ HHD of the aligned moving frames along the propagational direction on an atrium. (a) Irrotational vector with the potential $u$, (b) Incompressible vector with the potential $v$, and (c) harmonic vector with the potential $U$. }
\label {HHDAtrium}
\end{figure}

\begin{figure}[ht]
\centering
\subfloat[Irrotational] {\label{BunnyIrrot}  \includegraphics[width=4cm]{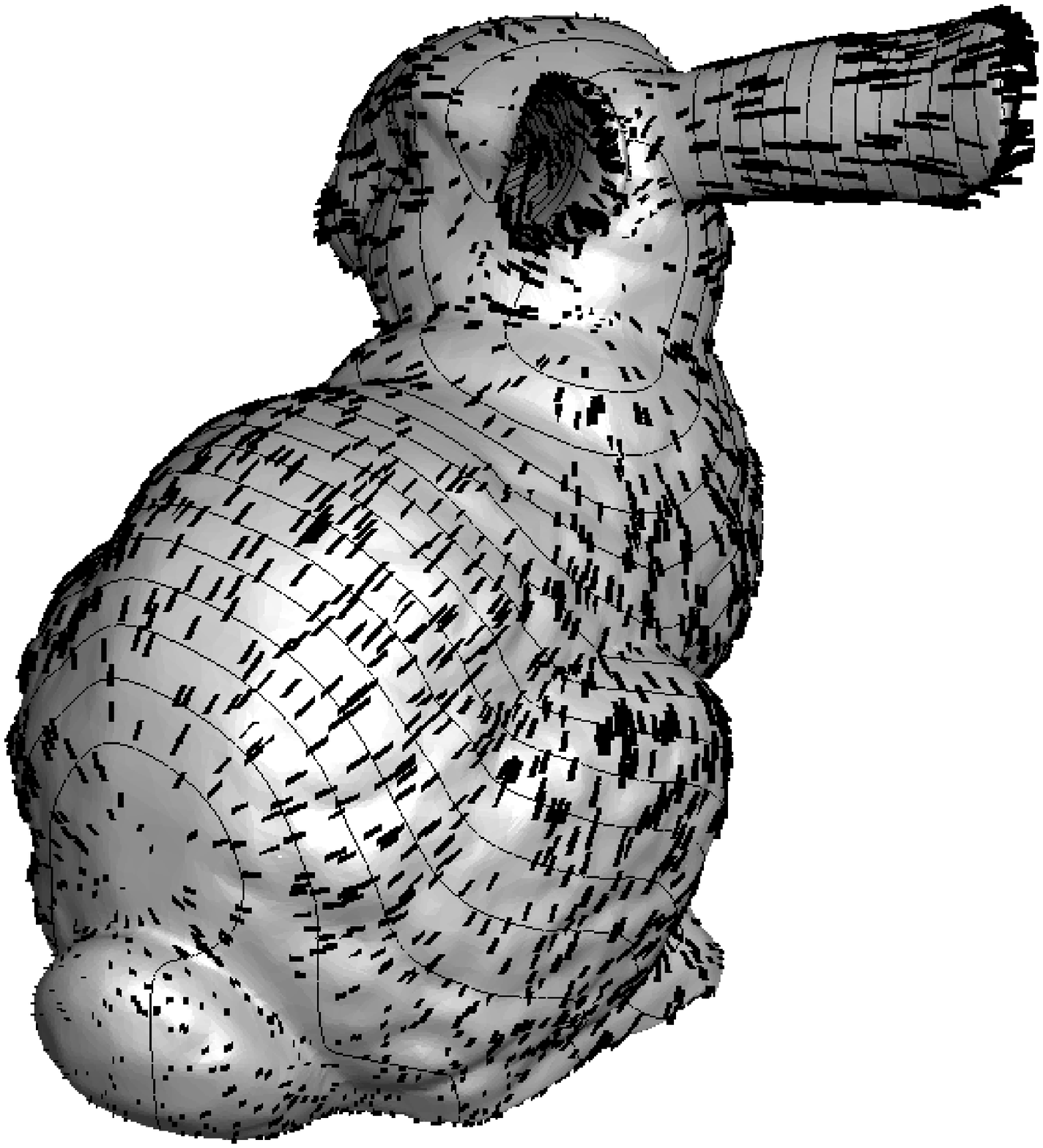}  }
\subfloat[Incompressible] {\label{BunnyIncomp}  \includegraphics[width=4cm]{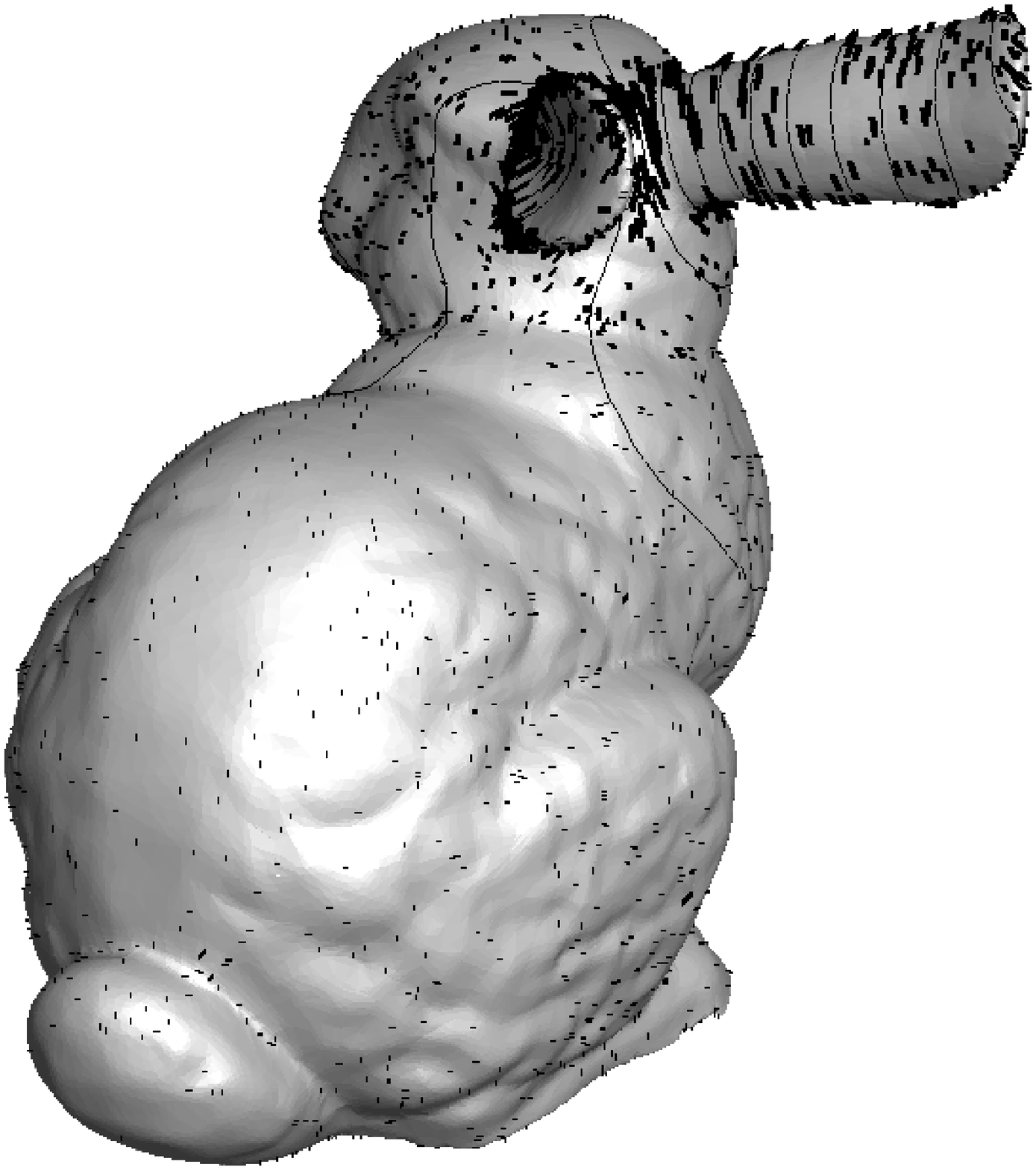}  }
\subfloat[Harmonic] {\label{Bunnyharmon}  \includegraphics[width=4cm]{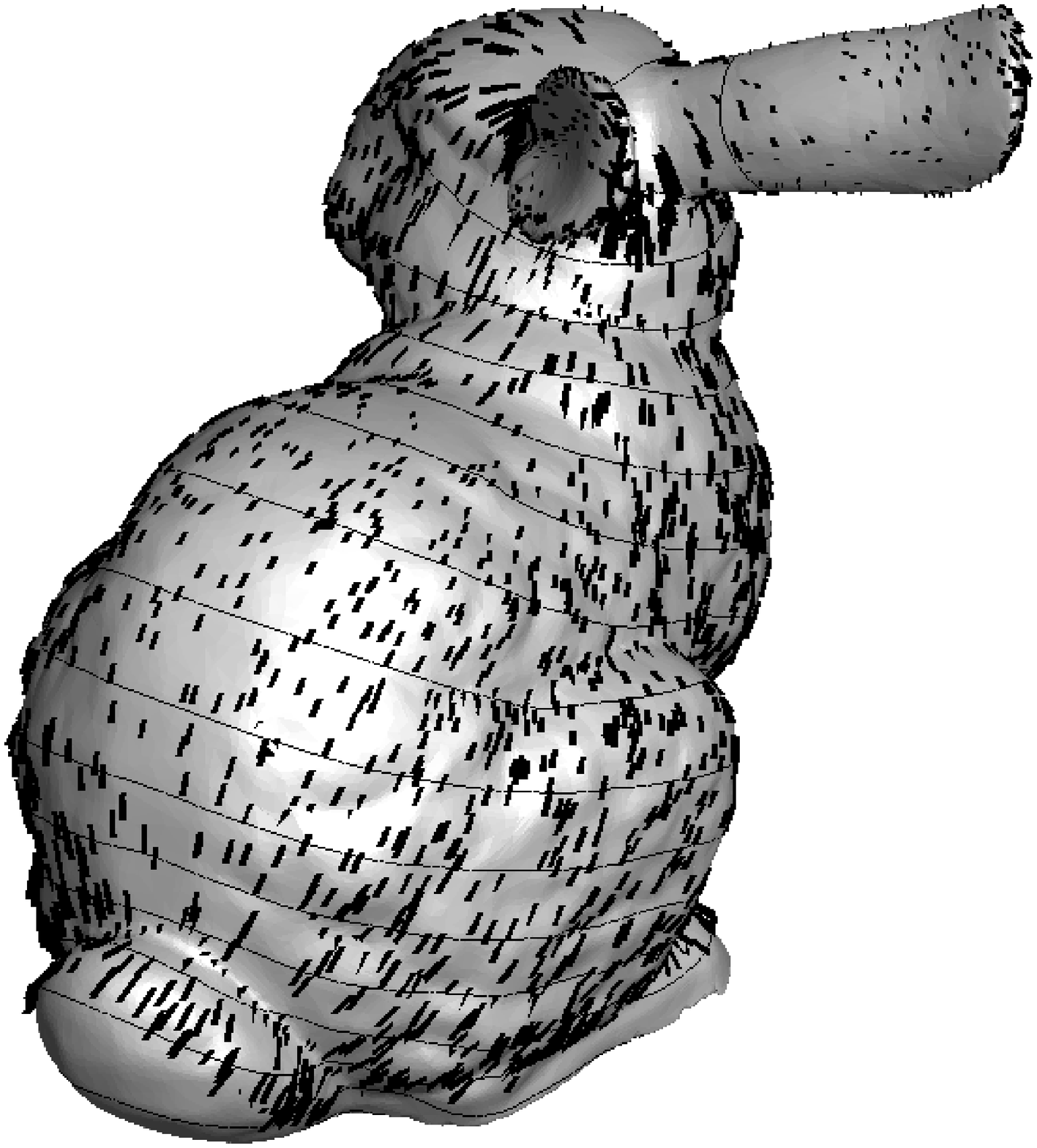}  }
\caption{ HHD of the aligned moving frames along the propagational direction on the Stanford Bunny. (a) Irrotational vector with the potential $u$, (b) Incompressible vector with the potential $v$, and (c) harmonic vector with the potential $U$.}
\label {HHDBunny}
\end{figure}

Two examples are used to demonstrate the proposed scheme even for a complexly-curved surface: the first is the surface model of the human atrium, and the second is the Stanford bunny. The initial vector is obtained by propagating a diffusion-reaction type wave from a point and by aligning moving frames along the gradient of the action potential \cite{MMFCNT}. Fig. \ref{HHDAtrium} and Fig. \ref{HHDBunny} present the HHD of the obtained vector into three components for the atrium and bunny, respectively. For the atrium, $\| \nabla \times \nabla u \| = 4.71e$-$13$, $\| \nabla \cdot  J \nabla v \| = 1.91e$-$12$ by two-dimensional discontinuous Helmsolver with moving frames. For the bunny, $\| \nabla \times \nabla u \| = 4.52e$-$11$, $\| \nabla \cdot  J \nabla v \| = 1.46e$-$09$ by two-dimensional continuous Helmsolver with moving frames. The magnitude of vector Laplacian is not negligible in some region of the domains, especially in the boundaries of the atrium and non-smooth junctions of the bunny, even though it still yields the smooth harmonic potential $U$. This problem could be the future work related to the development of the HHD in the context of Galerkin methods.

\section*{Acknowledgements}
This research was supported by the Basic Science Research Program through the National Research Foundation of Korea (NRF) and funded by the Ministry of Education, Science and Technology (No. 2016R1D1A1A02937255).

\bibliographystyle{elsarticle-num}
\bibliography{MMFCovariant}

\end{document}